\documentclass[11pt]{amsart}
\usepackage{amsfonts,amsmath,amssymb}

\newcommand{\bpf}[1][Proof]{{\noindent {\sc #1: }}}
\newcommand{\epf}{{{\hfill $\Box$ \smallskip}}}
\newtheorem{lemma}{Lemma}
\newtheorem{theorem}{Theorem}
\newtheorem{corollary}{Corollary}
\newtheorem{remark}{Remark}

\newcommand{\Pp}{\mathsf{P}}
\newcommand{\R}{\mathbb{R}}
\newcommand{\N}{\mathbb{N}}

\newcommand{\sgn}{\mathop{\rm sgn}}

\newcommand{\eps}{\varepsilon}

\newcommand{\Law}{\mathop{\mathsf{Law}}}

\title{Gumbel distribution in exit problems}
\author{Yuri Bakhtin}
\address{686 Cherry Street, School of Mathematics, Georgia Institute of Technology, Atlanta, GA, 30332-0160}
\email{bakhtin@math.gatech.edu}

\begin{document}

\maketitle

\begin{abstract} 
We explain the connection between the Gumbel limit for diffusion exit times and the theory of
extreme values.
\end{abstract}

\section{Introduction}
The Gumbel distribution $\Lambda$ is one of the max-stable distributions. Its
distribution function is
\begin{equation*}
\Lambda(x)=e^{-e^{-x}},\quad x\in\R,
\end{equation*}
and its density is given by
\begin{equation}
\lambda(x)=e^{-x-e^{-x}},\quad x\in\R.
\label{eq:Gumbel_density}
\end{equation}

The Gumbel distribution is mostly well-known as a limit law in the extreme value theory, see, e.g.,
\cite[Chapter 1]{deHaan:MR2234156}. Surprisingly, it also
appears as a limiting distribution (in the limit of vanishing noise)
for normalized exit times of diffusions conditioned on unlikely exit locations. The first result of
this kind was obtained in~\cite{Day:MR1110156} where exits through a characteristic repelling boundary were
considered. In a recent paper~\cite{ALEA} (see also references therein for related results) a similar
result was obtained for the case where the prescribed exit location for the diffusion is separated
from the starting point by a potential wall. Let us state the main result from~\cite{ALEA}.

 Let $X_{\eps}$ be a strong solution of the following one-dimensional 
stochastic differential equation
driven by a drift vector field $b$ and a Wiener process~$W$:
\begin{align}
\label{eq:basic-sde}
dX_\eps(t)&=b(X_\eps(t))dt+\eps dW(t),\\
X_\eps(0)&=x_0,
\end{align}
where $\eps>0$ is a constant diffusion coefficient and $x_0$ is a nonrandom starting point. 
This diffusion process is considered on a segment $[A,B]$ containing $0$ and $x_0$, until
the first exit time \[\tau_\eps=\inf\bigl\{t\ge 0: X_\eps(t)\in\{A,B\}\bigr\}.\]
It is assumed that the drift $b$ is smooth on $[A,B]$ and satisfies $b(0)=0$, $b'(0)>0$,
$b(x)<0$ for $x\in[A,0)$, and $b(x)>0$ for $x\in(0,B]$.  

Let us assume that $x_0<0$. Then the event $C_\eps=\{X_\eps(\tau_\eps)=B\}$ is unlikely for small $\eps$ 
because on $C_\eps$ the process $X_\eps$ has to travel against
the drift. The probability of $C_\eps$ decays exponentially in $\eps^{-2}$ as follows from the celebrated
Freidlin--Wentzell theory of large deviations for small white noise perturbations of dynamical systems, see~\cite{FW:MR722136}.
Nevertheless, one can study the behavior of~$\tau_\eps$ conditioned on this unlikely event.

\begin{theorem}[\cite{ALEA}] Under the conditions given above, there are constants $c_1,c_2,c_3$ such that
as $\eps\to0$, 
\[
\Law\left[\tau_\eps-c_1\ln\frac{1}{\eps}\ \Bigr|\ C_\eps\right]\ \Rightarrow\   \Law\left[c_2Z+c_3\right],
\]
where $Z$ is a 	Gumbel random variable and ``$\Rightarrow$'' denotes weak convergence
of probability measures.
\end{theorem}

Although it is clear that the behavior of both extreme values and exit times depends 
on the tails of the random variables involved, the precise reason why the Gumbel
distribution appears as a limit in both cases has not been known. 
The goal of this note is to explain the connection which is provided by the theory
of residual life times developed in the classical work~\cite{residual:MR0359049}. 

We stress that the computations in the present note are not new and go back at least to~\cite{Day:MR1110156}. 
What is new is the realization that these computations can be viewed as a specific case of those
in the heart of residual life times theory. In fact, this link between exit times and residual life times is absent
in the literature known to the author (and MathSciNet).

 In Section~\ref{sec:linear} we consider a simple model problem that 
simultaneously captures
the essence of the mechanism of generating 
asymptotically Gumbel exit times, and allows for elementary computations. 
In Section~\ref{sec:extreme-and-residual} we explain the necessary facts on extreme values and residual
life times. In Section~\ref{sec:log-transform}, we revisit the computations of Section~\ref{sec:linear}
and explain them
from the point of view of the theory of residual life times.

\medskip

{\bf Acknowledgements}. The author is grateful to Jon Wellner for pointing to~\cite{residual:MR0359049}. The author thanks
Gautam Goel for stimulating discussions and computer simulations. The partial support from NSF through CAREER Award DMS-0742424 
is gratefully acknowledged.

\section{A simple exit problem}\label{sec:linear}

Here we consider the following model problem. Let us suppose that the drift
is linear, i.e., $b(x)\equiv \beta x$ for an arbitrary
constant $\beta>0$ and, the constants $A$ and $B$ are equal to $-1$ and~$1$, respectively.
Instead of considering a fixed value of $x_0\in(-1,0)$, we will assume that $x_0=-\eps a$
for some $a>0$, and then sequentially take the limits $\eps\to0$ and $a\to+\infty$. 
This model is different from the one considered in~\cite{ALEA},\cite{Day:MR1110156},\cite{Day:MR1175267},\cite{Day:MR1376805},  and our results are not equivalent to the results therein. However, we choose this model and the limiting procedure 
since this is the simplest setting that suits our goal of demonstrating a natural connection between exit problems and
extreme values theory.

Let us first study the limiting behavior of $\tau_\eps$ as $\eps\to0$.
\begin{theorem}
\label{thm:eps-to-0} Let us fix $a>0$ and introduce $r=a\sqrt{2\beta}$.
As $\eps\to0$, 
\[
\Law\left[\tau_\eps-\frac{1}{\beta}\ln\frac{1}{\eps}\ \Bigr|\ X_\eps(\tau_\eps)=1\right]\Rightarrow
\Law\left[-\frac{1}{\beta}\ln(N-r)+\frac{1}{2\beta}\ln(2\beta)\ \Bigr|\ N>r\right],
\]
where $N$ is a standard Gaussian random variable.

\end{theorem}
\bpf In the case of linear drift, the solution $X_\eps$ can be represented by the following Duhamel principle:
\begin{align}
\notag
X_\eps(t)&=e^{\beta t}\left(x_0+\eps \int_0^te^{-\beta s}dW(s)\right)\\
&=e^{\beta t}\left(-\eps a+\eps \int_0^te^{-\beta s}dW(s)\right)
\label{eq:Duhamel}\\
\notag
&=\eps e^{\beta t}\left(-a+\int_0^te^{-\beta s}dW(s)\right).
\end{align}
conditioned on the exit through the right end of $[-1,1]$, i.e., on $\{X(\tau_\eps)=1\}$.

Let us introduce
\[I_t=\int_0^te^{-\beta s}dW(s),\quad t\in[0,+\infty],\]
and $I^*=\sup_{t\ge0} |I_t|$. Representing $I_t$ as a
time changed Wiener process, we see that
$I^*$ is a.s.-finite. Since $|X_\eps(\tau_\eps)|=1$, equation~\eqref{eq:Duhamel}
implies
\[
\tau_\eps\ge \frac{1}{\beta}\ln\frac{1}{\eps (I^*+a)},
\]
and we can conclude that $\tau_\eps\to+\infty$ a.s.
Combining this with~\eqref{eq:Duhamel} and continuity of $I_t$ on $[0,+\infty]$, we see that
\[
\frac{X_\eps(\tau_\eps)}{\eps e^{\beta \tau_\eps}}
=-a+I_{\tau_\eps}\stackrel{\mathrm{a.s.}}{\longrightarrow}-a+I_\infty.
\]
Recalling that $|X_\eps(\tau_\eps)|=1$, we conclude that
\[
\tau_\eps-\frac{1}{\beta}\ln\frac{1}{\eps}\stackrel{\mathrm{a.s.}}{\longrightarrow}-\frac{1}{\beta}\ln|-a+I_\infty|.
\]
Also,
\[
X_\eps(\tau_\eps)\stackrel{\mathrm{a.s.}}{\longrightarrow}\sgn(-a+I_\infty).
\]
Now the Theorem follows since $I_\infty$ is a centered Gaussian random variable with variance $1/(2\beta)$,
i.e., $I_\infty=N/\sqrt{2\beta}$ for a standard Gaussian random variable $N$. \epf

\medskip

The limiting distribution provided for the first part of the limiting  procedure ($\eps\to 0$) by Theorem~\ref{thm:eps-to-0} 
is a linear transformation of the conditional distribution 
$\Law\left[-\ln(N-r)\ \bigr|\ N>r\right]$, where $r=a\sqrt{2\beta}$.
In fact, this conditional distribution can be computed explicitly. For example, its density is given
by
\begin{equation}
\label{eq:p_r}
p_r(x)=\frac{\frac{1}{\sqrt{2\pi}}e^{-x-(e^{-x}+r)^2/2}}{1-G(r)},
\end{equation}
where $G$ is the standard Gaussian distribution function.

Similar arguments were used to study exits from neighborhoods of unstable equilibria in \cite{Day:MR1110156,
Day:MR1175267,Day:MR1376805,Berglund-Gentz:MR2039489, Bakhtin:MR2411523,Bakhtin:MR2731621,Bakhtin:MR2800902,Bakhtin:MR2935120,
Bakhtin:MR2983392, BG_periodic2}. 

\medskip

The second part of the analysis is letting $a\to\infty$ or, equivalently,
$r\to\infty$. Our next statement claims that as $r\to\infty$, under a proper normalization
(shift by~$\ln r$), the limit of the distributions described in Theorem~\ref{thm:eps-to-0}
is the Gumbel distribution. 

\begin{theorem}\label{thm:density_limit} Let $\lambda$ and $p_r$ be defined by~\eqref{eq:Gumbel_density} and~\eqref{eq:p_r}. Then 
\[
\lim_{r\to\infty} p_r(x-\ln r)= \lambda(x),\quad x\in\R.
\]
\end{theorem}
\bpf For any $x\in\R$,
\begin{align*} 
p_r(x-\ln r)&=\frac{\frac{1}{\sqrt{2\pi}}e^{-x-\ln r-(e^{-x-\ln r}+r)^2/2}}{1-G(r)}\\
&=\frac{\frac{1}{\sqrt{2\pi}r}e^{-r^2/2}}{1-G(r)} e^{-x- e^{-2x}/(2r^2)-e^{-x}}\to\lambda(x),\quad r\to\infty,
\end{align*}
where we used the fact that $1-G(r)\sim \frac{1}{\sqrt{2\pi}r}e^{-r^2/2}$.
\epf

\begin{corollary}\label{cor:log-conditioned-Gauss} As $r\to\infty$,
\[
\Law\left[-\ln(N-r)-\ln r\ \bigr|\ N>r\right]\Rightarrow \Lambda.
\]
\end{corollary}
\bpf  This claim is a direct consequence of Theorem~\ref{thm:density_limit} and Scheff\'e's theorem (see~\cite{Scheffe:MR0021585} or \cite[Appendix II]{Billingsley:MR0233396}) stating that pointwise convergence
of densities implies convergence of distributions in total variation and hence weak convergence. One can also give a direct proof:
\begin{align}\label{eq:proof_of_conv_to_gumbel}
\frac{\Pp\{-\ln(N-r)-\ln r<x\}}{\Pp\{N>r\}}&=\frac{\Pp\{N>r+e^{-x}/r\}}{\Pp\{N>r\}}\\
&\sim\frac{
\frac{1}{\sqrt{2\pi}(r+e^{-x}/r)}e^{-(r+e^{-x}/r)^2/2}
}{
\frac{1}{\sqrt{2\pi}r}e^{-r^2/2}}\sim e^{-e^{-x}}
\notag
\end{align}
\epf

Although this convergence statement along with Theorem~\ref{thm:density_limit} is a result of a straightforward
computation that in some form appeared for the first time in the diffusion context in~\cite{Day:MR1110156}, 
we still need to explain the connection to the theory of extreme values and residual life times, and we proceed to recall 
some related basic facts.

\section{Gumbel distribution in extreme value theory and residual life time theory}
\label{sec:extreme-and-residual}
Let us start with some results from \cite{Gnedenko:MR0008655} where basins of attraction of $\Lambda$ and
other max-stable distributions were studied. If random variables $X_1,X_2,\ldots$ are i.i.d.\ with common distribution function $F$,
then
\[
\Pp\{\max(X_1,\ldots,X_n)\le x\}=F^n(x),\quad x\in\R, n\in\N.
\]
The following is a version of
Lemma~3 in \cite{Gnedenko:MR0008655}:
\begin{lemma}\label{lem:Gnedenko-criterion}
Let
 $F$ be a distribution function,
$\Phi$ a continuous distribution function, 
and let $(a_n), (b_n)$ be number sequences, $a_n>0$ for all $n$. Then
\begin{equation}
\lim_{n\to\infty}F^n(a_nx+b_n)= \Phi(x),\quad  x\in\R,
\label{eq:G33}
\end{equation}
if and only if
\begin{equation}
\Phi(x)\ne 0\quad \Rightarrow\quad \lim_{n\to\infty} n[1-F(a_n x+b_n)] = -\ln \Phi(x).
\label{eq:G34}
\end{equation}
\end{lemma}

\bpf Let us recall the proof of this statement from \cite{Gnedenko:MR0008655}. Condition~\eqref{eq:G33}
is equivalent to 
\begin{equation}
\Phi(x)\ne 0\quad\Rightarrow\quad\lim_{n\to\infty} n\ln F(a_nx+b_n)=\ln \Phi(x).
\label{eq:G36}
\end{equation}
To deal with the l.h.s.\ of~\eqref{eq:G36} we notice that \eqref{eq:G33} implies
\begin{equation}
\Phi(x)\ne 0\quad\Rightarrow\quad\lim_{n\to\infty} F(a_nx+b_n)
=1,
\label{eq:G35}
\end{equation}
which in turn implies
\begin{equation}
\label{eq:G37}
\Phi(x)\ne 0\quad\Rightarrow\quad \ln F(a_nx+b_n)= -(1-F(a_nx+b_n))(1+o(1)),\quad n\to\infty.
\end{equation}
Now,~\eqref{eq:G34} follows from~\eqref{eq:G36} and~\eqref{eq:G37}.

If~\eqref{eq:G34} holds, then so does~\eqref{eq:G35} and, therefore,~\eqref{eq:G37}, which implies~\eqref{eq:G36} which
is equivalent
to~\eqref{eq:G33}.
\epf


\medskip

If~\eqref{eq:G33} holds for some sequences $(a_n)$ and $(b_n)$ with $a_n>0$, then we will write $F\in D(\Phi)$.
%

Let $G$ be the distribution function of the standard Gaussian distribution. The following simple
result shows that $G\in D(\Lambda)$.
\begin{lemma}\label{lem:Gauss-tail}
Let $b_n$ satisfy
$1-G(b_n)=n^{-1}$
and let $a_n=b_n^{-1}$. Then
\[
\lim_{n\to\infty}G^n(a_n x+b_n)=\Lambda(x),\quad x\in\R.
\]
\end{lemma}
\bpf Noticing that as $n\to\infty$, one has $b_n\to\infty$, $a_n\to 0$, $a_nx+b_n\to\infty$,  we can write 
\begin{align*}
n(1-G(a_nx+b_n))&\sim n\frac{1}{\sqrt{2\pi} (a_nx+b_n)}e^{-(a_nx+b_n)^2/2}\\
&\sim n\frac{1}{\sqrt{2\pi} b_n}e^{-b_n^2/2}e^{-b_na_nx}e^{-a_n^2x^2/2}\\
&\sim n(1-G(b_n))e^{-x}\sim e^{-x},
\end{align*}
i.e.,
\[
\lim_{n\to\infty}n(1-G(a_nx+b_n)) =-\ln \Lambda(x),\quad x\in\R,
\]
and the result follows from Lemma~\ref{lem:Gnedenko-criterion}.\epf

\bigskip

Let us now turn to residual life times. For a random variable $X$ with distribution function $F(x)=\Pp\{X\le x\}$ and tail 
\begin{equation}
\label{eq:tail}
R(x)=\Pp\{X>x\}=1-F(x),
\end{equation}
the distribution tail of residual life time after time $t$ is
\begin{equation}
\label{eq:residual_time}
R_r(x)=\Pp\{X-r>x| X>r\}=\frac{R(r+x)}{R(t)}.
\end{equation}
In~\cite{residual:MR0359049}, a theory of scaling limits
for residual life times was developed and connections with
the theory of scaling limits of extreme values were established.
The following is a version of
Theorem~3 from~\cite{residual:MR0359049}.
\begin{theorem}[\cite{residual:MR0359049}]
\label{thm:extreme-vs-residual} Let $F$ be a distribution function,  let $R$ and $R_r$ be defined by~\eqref{eq:tail} and \eqref{eq:residual_time}, and let $(a_n), (b_n)$ be number sequences satisfying $a_n>0$ for all $n$.
Then the following two conditions are equivalent:
\begin{enumerate}
\item \label{it:max-limit} For all $x\in\R$, $F(x)<1$ and $\lim_{n\to\infty}F^n(a_nx+b_n)=\Lambda(x)$.
\item \label{it:residual-limit} There are functions $a(r)>0$ and $b(r)$ such that
\begin{equation*}
\lim_{n\to\infty} \frac{R(a(r)x+b(r))}{R(r)}= -\ln \Lambda(x)=e^{-x},\quad x\in\R.
\end{equation*}

\end{enumerate}
\end{theorem}
\bpf[Sketch of proof] Suppose condition~\ref{it:max-limit} is satisfied. Lemma~\ref{lem:Gnedenko-criterion} implies that
\begin{equation}
\lim_{n\to\infty} nR(a_nx +b_n) = -\ln \Lambda(x)=e^{-x},\quad x\in\R,
\label{eq:conv_to_ln_Lambda}
\end{equation}
and condition~\ref{it:residual-limit} follows with $a$ and $b$ satisfying $a(r)=a_n$ and $b(r)=b_n$ for all 
$r$ such that $(n+1)^{-1}\le R(r)<n^{-1}$. We omit the proof of the converse implication and only mention that it is also based on~\eqref{eq:conv_to_ln_Lambda}. 
\epf

\begin{remark}\label{rem:a(t)}\rm In fact, the function $b$ in condition~\ref{it:residual-limit} of the theorem can be chosen to be  $b(r)=r$
(see \cite{residual:MR0359049}, Corollary~2),
so that condition~\ref{it:residual-limit} can be rewritten as
\begin{equation}
\label{eq:b_eq_t}
\lim_{n\to\infty}\frac{R(r+a(r)x)}{R(r)}= -\ln \Lambda(x) = e^{-x},\quad x\in\R.
\end{equation}
\end{remark}
\section{Logarithmic transformation of residual life times}\label{sec:log-transform}

Let $X$ be a random variable with distribution function $F$ and tail function $R=1-F$.
Motivated by the results of Section~\ref{sec:linear}, we would like to study the limiting behavior of 
$\Law[-\ln (X-r)\ |\ X>r]$, so we
are interested in
\[
H_r(x)=\Pp\{-\ln(X-r)\le x| X>r\}.
\]
\begin{theorem}\label{thm:log-transform} Let $F\in D(\Lambda)$. Let $R$ be defined by~\eqref{eq:tail}. 
If~\eqref{eq:b_eq_t} holds for some function $a(\cdot)$, then
\[
\lim_{n\to\infty}H_r(x-\ln a(r))=\Lambda(x).
\]
\end{theorem}
\bpf We can rewrite $H_r(x)$ as
\[
H_r(x)=\Pp\{X-r> e^{-x}| X>r\}=\frac{R(r+e^{-x})}{R(r)}.
\]
Since
\[
H_t(x-\ln a(r))=\frac{R(r+a(r)e^{-x})}{R(r)},
\]
Theorem~\ref{thm:extreme-vs-residual}, Remark~\ref{rem:a(t)},
and the basic identity
\begin{equation}
\label{eq:invariance-of-Lambda}
-\ln\Lambda(e^{-x})=\Lambda(x),
\end{equation}
imply our claim.
\epf

\bigskip

Theorem~\ref{thm:log-transform} can be applied to the Gaussian distribution as the following lemma shows.

\begin{lemma}\label{lem:Gaussian-tail} If $F$ is the standard Gaussian distribution function $G$, then 
one can choose $a(t)=t^{-1}$ in~\eqref{eq:b_eq_t}.
\end{lemma}
\bpf The following computation for $R=1-G$ is essentially the same as~\eqref{eq:proof_of_conv_to_gumbel}
in the proof of Corollary~\ref{cor:log-conditioned-Gauss}: 
\begin{align*}
\frac{R(t+x/t)}{R(t)}\sim\frac{\frac{1}{\sqrt{2\pi}(t+x/t)}e^{-(t+x/t)^2/2}}{\frac{1}{\sqrt{2\pi}t}e^{-t^2/2}}
\sim e^{-x},\quad t\to\infty.
\end{align*} \epf

Now Corollary~\ref{cor:log-conditioned-Gauss} responsible for the Gumbel limit for exit times can be seen
as a direct consequence of Theorem~\ref{thm:log-transform} and Lemma~\ref{lem:Gaussian-tail}. This puts
the Gumbel limit of exit times in the context of the residual life time theory. Interestingly, this
connection of exit times to the residual times theory depends on a trivial but seemingly random
property~\eqref{eq:invariance-of-Lambda} of the Gumbel distribution.

\bibliographystyle{alpha}
\bibliography{gumbel}

\end{document}